\documentstyle{article}

\newcommand{\rto}{\rightarrow}

\newcommand{\dotcup}{+}

\newcommand{\EOP}{\Box}

\newtheorem{theorem}{Theorem}[section]
\newtheorem{lemma}[theorem]{Lemma}
\newtheorem{corollary}[theorem]{Corollary}
\newtheorem{proposition}[theorem]{Proposition}
\newtheorem{definition}[theorem]{Definition}
\newtheorem{conjecture}[theorem]{Conjecture}

\input epsf

\title{Tripartite Version of the Corr\'adi-Hajnal Theorem}
\author{Csaba Magyar and Ryan R. Martin\footnote{Carnegie Mellon
University, Pittsburgh, PA.  Email: {\bf
rymartin@andrew.cmu.edu}}~\footnote{This author partially supported by
DIMACS via NSF grant CCR 91-19999}}
\date{}

\begin{document}

\maketitle

\begin{abstract}
Let $G$ be a tripartite graph with $N$ vertices in each vertex class.
If each vertex is adjacent to at least $(2/3)N$ vertices in each of
the other classes, then either $G$ contains a subgraph that consists
of $N$ vertex-disjoint triangles or $G$ is a specific graph in which
each vertex is adjacent to exactly $(2/3)N$ vertices in each of the
other classes.
\end{abstract}

\section{Introduction}
\label{sINTRO}

A central question in extremal graph theory is the determination of
the minimum density of edges in a graph $G$ which guarantees a
monotone property $\cal P$.  If the property is the inclusion of a
fixed size subgraph $H$, the answer is given by the classic theorems
of Tur\'{a}n~\cite{T1} (when $H$ is a complete graph) and Erd\H{o}s
and Stone~\cite{ES1}.

However, in the case when a graph $G$ is required to contain a
spanning subgraph $H$; that is, $H$ has the same number of vertices as
$G$, an important parameter is a lower bound on the minimum degree
that guarantees $H$ is a subgraph of $G$.  Perhaps the most well-known
result of this type is a theorem of Dirac~\cite{D1} which asserts that
every $n$ vertex graph with minimum degree at least $\frac{n}{2}$
contains a Hamiltonian cycle.  Another theorem of this type is the
so-called Hajnal-Szemer\'{e}di theorem, with the case $k=3$ proven
first by Corr\'adi and Hajnal~\cite{CH1}.
\begin{theorem}[Hajnal-Szemer\'{e}di~\cite{HS1}]  Let $G$ be
graph on $n$ vertices with minimum degree $\frac{k-1}{k}n$.  If $k$
divides $n$, then $G$ has a subgraph that consists of $\frac{n}{k}$
vertex-disjoint cliques of size $k$.
\label{tHAJ-SZEM}
\end{theorem}

A tripartite graph is said to be balanced if it contains the same
number of vertices in each class.  Theorem~\ref{tTRI} is a tripartite
version of the Corr\'adi-Hajnal result.
\begin{theorem}
Let $G=(V_1,V_2,V_3;E)$ be a balanced tripartite graph on $3N$
vertices such that each vertex is adjacent to at least $(2/3)N$
vertices in each of the other classes.  If $N\geq N_0$ for some
absolute constant $N_0$, then $G$ has a subgraph consisting of $N$
disjoint triangles or $G=\Gamma_3(N/3)$ for $N/3$ an odd integer.
\label{tTRI}
\end{theorem}



The graph $\Gamma_3(N/3)$ is defined in Section~\ref{sDEFS}.  The
proof is in two parts.  Theorem~\ref{tTRIFUZZY} in
Chapter~\ref{sTRIFUZZY} states that if the degree condition is
relaxed, then all graphs, except a specific class, have the spanning
subgraph of disjoint triangles.  We will then show how find the
spanning subgraph for that excluded class of graphs by proving
Theorem~\ref{tTRIEXTREME} in Chapter~\ref{sTRIEXTREME}.  Assume that
$N$ is divisible by 3.  If not, Section~\ref{sNOTDIVIS} shows that the
case where $N$ is not divisible by $3$ comes as a corollary.

\subsection{The Regularity and Blow-up Lemmas}
Throughout this paper, we will try to keep much of the notation and
definitions in~\cite{KSS1}.  The symbol $\dotcup$ will sometimes be
used to denote the disjoint union of sets.  $V(G)$ and $E(G)$ denote
the vertex-set and edge-set of the graph $G$, respectively.  The
triple $(A,B;E)$ denotes a bipartite graph $G=(V,E)$, where
$V=A\dotcup B$ and $E\subset A\times B$.  $N(v)$ denotes the set of
neighbors of $v\in V$.  For $U\subset V\setminus\{v\}$, $N_U(v)$
denotes the set of neighbors of $v$ intersected with $U$.  The {\bf
degree of $v$} is $\deg(v)=|N(v)|$.  The {\bf degree of $v$ in $U$} is
$\deg_U(v)=\left|N_U(v)\right|$.  If $H$ is a subgraph of $G$, then we
relax notation so that $\deg_H(v)=\deg_{V(H)}(v)$.  For $U\subset
V$, $\left. G\right|_U$ denotes the graph $G$ induced by the vertices
$U$.

The graph $K_3$ is the complete graph on $3$ vertices, the ``{\bf
triangle}.''  We say edges and triangles are {\bf disjoint} if their
common vertex set is empty.  A balanced tripartite graph on $3N$
vertices is {\bf covered with triangles} if it contains a subgraph of
$N$ disjoint triangles.  The tripartite version of the
Corr\'adi-Hajnal result is Theorem~\ref{tTRI}.  When $A$ and $B$ are
subsets of $V(G)$, we define
\[ e(A,B)=\left|\left\{(x,y) : x\in A, y\in B, \{x,y\}\in
   E(G)\right\}\right| . \] 
For nonempty $A$ and $B$,
\[ d(A,B)=\frac{e(A,B)}{|A||B|} \]
is the {\bf density} of the subgraph of edges that contain one
endpoint in $A$ and one in $B$.

\begin{definition} The bipartite graph $G=(A,B,E)$ is {\bf
$\epsilon$-regular} if 
\[ X\subset A, Y\subset B, |X|>\epsilon|A|, |Y|>\epsilon|B| \]
imply $|d(X,Y)-d(A,B)|<\epsilon$, otherwise we say $G$ is {\bf
$\epsilon$-irregular}.
\label{def1}
\end{definition}

We will also need a stronger version.
\begin{definition} $G=(A,B,E)$ is {\bf
$(\epsilon,\delta)$-super-regular} if 
\[ X\subset A, Y\subset B, |X|>\epsilon|A|, |Y|>\epsilon|B| \]
imply $d(X,Y)>\delta$ and
\[ \deg(a)>\delta|B|, \quad \forall a\in A\qquad\mbox{and}\qquad
\deg(b)>\delta|A|, \quad \forall b\in B . \]
\label{def2}
\end{definition}

One of our main tools will be the Regularity Lemma~\cite{Sz1}, but
more specifically, a corollary known as the Degree Form:
\begin{lemma}[Degree Form of the Regularity Lemma] 
For every positive $\epsilon$ there is an $M=M(\epsilon)$ such that if
$G=(V,E)$ is any graph and $d\in [0,1]$ is any real number, then there
is a partition of the vertex set $V$ into $\ell+1$ clusters
$V_0,V_1,\ldots,V_{\ell}$ and there is a subgraph $G'=(V,E')$ with the
following properties:
\begin{itemize}
\item $\ell\leq M$,
\item $|V_0|\leq\epsilon|V|$,
\item all clusters $V_i$, $i\geq 1$, are of the same size
$L\leq\left\lceil\epsilon|V|\right\rceil$,
\item $\deg_{G'}(v)>\deg_G(v)-(d+\epsilon)|V|$, $\forall v\in V$,
\item $\left. G'\right|_{V_i}=\emptyset$ ($V_i$ are independent in
$G'$),
\item all pairs $\left. G'\right|_{V_i\times V_j}$, $1\leq i<j\leq l$,
are $\epsilon$-regular, each with density either $0$ or exceeding
$d$.
\end{itemize}
\label{reglem}
\end{lemma}


The above definition is the traditional statement of the Degree Form.
In fact, we can guarantee that each cluster that is not $V_0$ has that
all vertices belong to the same vertex class.  The Degree Form is
derived from the original Regularity Lemma (see~\cite{KSS3}) which
shows that any partition can be refined so that it is in the form of
the Regularity Lemma.  The reduced graph $G_r$, has a vertex set
$V_1,\ldots,V_{\ell}$ with $V_i\sim V_j$ if and only if
$\left. G'\right|_{V_i\times V_j}$ is $\epsilon$-regular with density
exceeding $d$.

We will also make use of the so-called Blow-up Lemma.  The graph $H$
can {\bf be embedded into} graph $G$ if $G$ contains a subgraph
isomorphic to $H$.
\begin{lemma}[Blow-up Lemma~\cite{KSS2}]  Given a graph $R$ of order
$r$ and positive parameters $\delta$, $\Delta$, there exists an
$\epsilon>0$ such that the following holds:  Let $N$ be an arbitrary
positive integer, and let us replace the vertices of $R$ with pairwise
disjoint $N$-sets $V_1,V_2,\ldots,V_r$ (blowing up).  We construct two
graphs on the same vertex-set $V=\cup V_i$.  The graph $R(N)$ is
obtained by replacing all edges of $R$ with copies of the complete
bipartite graph $K_{N,N}$ and a sparser graph $G$ is constructed by
replacing the edges of $R$ with some $(\epsilon,\delta)$-super-regular
pairs.  If a graph $H$ with maximum degree $\Delta(H)\leq\Delta$ can
be embedded into $R(N)$, then it can be embedded into $G$.
\label{blowuplem}
\end{lemma}

\subsection{Further Definitions}
\label{sDEFS}
We will frequently refer to the well-known K\"{o}nig-Hall condition,
which states that if $G=(A,B;E)$ is a bipartite graph, then there is a
matching in $G$ that involves all the vertices of $A$ unless there
exists an $X\subseteq A$ such that, $|N(X)|<|X|$.  Specifically, we
often use the immediate corollary that if $|A|=|B|$ and each vertex in
$A$ has degree at least $|B|/2$ and each vertex in $B$ has degree at
least $|A|/2$, then $G$ must have a perfect matching.

With $G$ a $k$-partite graph, $V(G)=V_1+\cdots+V_k$, each $V_i$ being
a partition, we refer to each $V_i$ as a {\bf vertex class}.  We refer
to the graph defined by the Regularity Lemma, denoted $G_r$, as the
{\bf reduced graph} of $G$.  $G$ itself is the {\bf real graph}.  Any
triangle in $G_r$ or in a similar reduced graph is referred to as a
{\bf super-triangle}.  A triangle in $G$ is often called a {\bf real
triangle}~to avoid confusion.

The notation $a\ll b$ means that the constant $a$ is small enough
relative to $b$.  This has become standard notation in these kinds of
proofs.  A set is of size {\bf $\gamma$-approximately $M$} if its size
is $(1\pm \gamma)M$.  Let us also define two classes of graphs.  The
first is $\Theta_{m\times n}$.  The vertices of $\Theta_{m\times n}$
are $\{h_{i,j} : i=1,\ldots,m; j=1,\ldots,n\}$ and $h_{i,j}\sim
h_{i',j'}$ iff $i\neq i'$ and $j\neq j'$.  Note that $\Theta_{3\times
2}$ contains no triangle.  The second graph is the graph $\Gamma_k$.
The vertices are $\{h_{i,j} : i=1,\ldots,k; j=1,\ldots,k\}$ and the
adjacency rules are as follows: $h_{i,j}\sim h_{i',j'}$ if $i\neq i'$
and $j\neq j'$ and either $j$ or $j'$ is in $\{1,\ldots,k-2\}$.  Also,
$h_{i,k-1}\sim h_{i',k-1}$ and $h_{i,k}\sim h_{i',k}$ for $i\neq i'$.
No other edges exist.  If $k$ is even, then $\Gamma_k$ can be covered
by $K_k$'s, but it cannot if $k$ is odd.

For a graph $G$, define $G(t)$ to be the graph formed by replacing
each vertex with a cluster of $t$ vertices and each edge with the
complete bipartite graph $K_{t,t}$.  For $\epsilon\geq 0$ and
$\Delta\geq 0$, a graph $H$ is {\bf $(\epsilon,\Delta)$-approximately
$G(t)$} if each vertex of $G$ is replaced with a cluster of size
$\epsilon$-approximately $t$ and each non-edge is replaced by a
bipartite graph of density at most $\Delta$.  For brevity, we will say
a graph is {\bf $\Delta$-approximately $G(t)$} if it is
$(0,\Delta)$-approximately $G(t)$.  Note that if $\Delta<\Delta'$ and
$\epsilon\ll\Delta'-\Delta$, then (if we are allowed to add or
subtract vertices to guarantee that clusters are the same size) a
graph that is $(\epsilon,\Delta)$-approximately $G(t)$ is also
$\Delta'$-approximately $G(t)$.

\subsection{An Easy Result}

Let $G$ be a balanced tripartite graph on $3M$ vertices such that each
vertex in $G$ is adjacent to at least $(3/4)M$ vertices in each of the
other classes.  Proposition~\ref{pEASYCASE} shows that this graph can
be covered with triangles.  Proposition~\ref{pEASYCASE} is used
repeatedly in Section~\ref{sTRIEXTREME}.
\begin{proposition}
Let $G=(V_1,V_2,V_3;E)$ be a balanced tripartite graph on $3M$
vertices such that each vertex is adjacent to at least $(3/4)M$
vertices in each of the other classes.  Then, we can cover $G$ with
$M$ vertex-disjoint triangles.
\label{pEASYCASE}
\end{proposition}

\noindent{\bf Proof.} Let $H$ be the graph induced by $(V_2,V_3)$.
Each vertex in $H$ is adjacent to at least $(3/4)M>(1/2)M$ vertices in
each of the other classes.  Therefore, $H$ can be covered by $M$
disjoint edges.  Each of these edges is adjacent to at least
$\left(1-2\times\frac{1}{4}\right)M=M/2$ vertices in $V_1$ and each
vertex in $V_1$ is adjacent to at least $M/2$ of the disjoint edges.
So, by K\"onig-Hall, there exists a 1-factor between $V_1$ and the $M$
disjoint edges -- giving us our $M$ disjoint triangles. $\EOP$

\subsection{A Useful Proposition}
\label{sTHETA}
Proposition~\ref{pTHETA} is quite valuable and is used in both
Section~\ref{sTRIFUZZY} and Section~\ref{sTRIEXTREME}.
\begin{proposition}
For a $\Delta$ small enough, there exists $\epsilon>0$ such that if
$H$ is a tripartite graph with at least $2\left(1-\epsilon\right)t$
vertices in each vertex class and each vertex is nonadjacent to at
most $\left(1+\epsilon\right)t$ vertices in each of the other classes.
Furthermore, let $H$ contain no triangles.  Then, each vertex class is
of size at most $2\left(1+\epsilon\right)t$ and $H$ is
$(\epsilon,\Delta)$-approximately $\Theta_{3\times 2}(t)$.
\label{pTHETA}
\end{proposition}

\noindent{\bf Proof.} Let $\epsilon\ll\delta\ll\delta'\ll\Delta$.
First we bound the sizes of the $V_i$.  Choose vertices $v_1$ and
$v_2$ from $V(G)\setminus V_i$ such that they form an edge.  These
vertices can have no common neighbor, giving that $|V_i|\leq
2(1+\epsilon)t$.

Now choose $w\in V_3$.  Let $N(w)\cap V_i$ be written as $A_{i,1}$,
for $i=1,2$, such that each vertex in $A_{i,1}$ is adjacent to no
vertices in $A_{3-i,1}$.  Furthermore, define $A_{3,1}$ to be those
vertices in $V_3$ that are adjacent to less than $\delta t$ vertices
in each of $A_{1,1}$ and $A_{2,1}$.  The set $A_{3,1}$ cannot be of
size larger than $(1+\epsilon) t$.  If it were, then there exists an
edge in $(A_{2,1},A_{3,1})$.  By the degree condition, if $\delta$ is
small enough, this edge must have a common neighbor in $V_1\setminus
A_{1,1}$.

For all $i\in [3]$, remove vertices (if necessary) from the sets
$A_{i,1}$ to create $A_{i,1}'$ so that each vertex in $A_{i,1}'$ is
adjacent to less than $\delta t$ vertices in each $A_{i',1}$ for
$i'\neq i$.  By the same arguments given before, $|A_{i,1}'|\leq
(1+\epsilon) t$, for $i=1,2,3$.  As a result, each vertex in
$A_{i,1}'$ is adjacent to less than $\delta' t$ vertices in each
$A_{i',1}'$, for $i'\neq i$.  Let $A_{i,2}'=V_i\setminus A_{i,1}'$ for
$i=1,2,3$.

We now want to show that each pair of the form
$\left(A_{i,2}',A_{i',2}'\right)$ is sparse.  Let $v\in A_{1,2}'$.  If
$N(v)\cap A_{2,2}'\neq\emptyset$, then $\left|A_{3,1}\setminus
N(v)\right|\leq\delta t$ which implies $\left|A_{2,1}\setminus
N(v)\right|\leq\delta t$.  As a result, $\left|N(v)\cup
A_{2,1}\right|,\left|N(v)\cup A_{3,1}\right|\leq(1+\epsilon+\delta)t$,
implying $\left|N(v)\cap A_{2,2}'\right|,\left|N(v)\cap
A_{3,2}'\right|\leq\delta' t$.  Similar results occur for $w\in
A_{2,2}'\cup A_{3,2}'$.  Once again, it must be true that each
$\left|A_{i,2}'\right|\leq(1+\epsilon)t$.  

Note that each set $A_{i,j}'$ is of size at least $(1-3\epsilon)t$
because the others are of size at most $(1+\epsilon)t$.  Therefore,
vertices can be moved from the sets larger than $(1-\epsilon)t$ to the
smaller sets to create sets $A_{i,j}''$ of size in
$\left((1-\epsilon)t,(1+\epsilon)t\right)$ with pairwise density at
most $\Delta$. $\EOP$

\newcommand{\caa}{3}
\newcommand{\caap}{4}
\newcommand{\threetimescaa}{9}
\newcommand{\cbb}{15}
\newcommand{\caapcbb}{18}

\section{The Fuzzy Tripartite Theorem}
\label{sTRIFUZZY}
\subsection{Statement of the Theorem}

Theorem~\ref{tTRIFUZZY} allows us, with an exceptional case, to cover
$G$ with triangles, even if the minimum degree is a bit less than
$(2/3)N$.

\begin{theorem}
Given $\epsilon\ll\Delta\ll 1$, let $G=(V_1,V_2,V_3;E)$ be a balanced
tripartite graph on $3N$ vertices such that each vertex is adjacent to
at least $(2/3-\epsilon)N$ vertices in each of the other classes.
Then, if $N$ is large enough, either $G$ can be covered with
triangles, or $G$ has three sets of size $N/3$, each in a different
vertex class, with pairwise density at most $\Delta$.
\label{tTRIFUZZY}
\end{theorem}

\subsection{Proof of the Theorem}
As usual, there is a sequence of constants:
\[ \epsilon\ll\epsilon_1\ll\epsilon_5\ll
   \epsilon_3\ll\alpha\ll\delta_4\ll\delta_3\ll d_3\ll d_1\ll
   \epsilon_2\ll\Delta_0\ll\Delta \]

Begin with $G=(V_1,V_2,V_3;E)$, a balanced tripartite graph on $3N$
vertices with each vertex adjacent to at least $(2/3-\epsilon)N$
vertices in each of the other classes.  Define the extreme case to be
the case where $G$ has three sets of size $N/3$ with pairwise density
at most $\Delta$.  Apply the Degree Form of the Regularity Lemma
(Lemma~\ref{reglem}), with $d_1$ and $\epsilon_1$, to partition each
of the vertex classes into $\ell+\caap$ clusters.  Let us define
$G_r'$ to be the reduced graph defined by the Lemma.  It may be
necessary to place clusters into the exceptional sets (the sets of
vertices in each vertex class that make up the $V_0$ in
Lemma~\ref{reglem}) to ensure that $\ell$ is divisible by 3.  It is
important to observe that in the proof, the exceptional sets will
increase in size, but will always remain of size $O(\epsilon_1)N$.

For $i=1,2,3$, there exist
$V_i=V_i^{(0)}+V_i^{(1)}+\cdots+V_i^{(\ell+\caa)}$ and
$\left|V_i^{(j)}\right|=L\leq\lceil\epsilon_1 N\rceil$, $\forall i$,
$\forall j\geq 1$.  The reduced graph $G_r'$ has the condition that
every cluster is adjacent to at least $(2/3-\epsilon_2)(\ell+\caa)$
clusters in each of the other vertex classes.  Apply
Lemma~\ref{lALMOSTCOVER} repeatedly to $G_r'$ with $M=\ell+\caa$ to
get a decomposition of $G_r'$ into $\ell$ vertex-disjoint triangles.
If this is not possible, then Lemma~\ref{lALMOSTCOVER} and
Proposition~\ref{pEXTREME} imply that $G$ is in the extreme case.
\begin{lemma} [Almost-covering Lemma]
Let $\epsilon'\ll\Delta_0\ll 1$, and let $G=(V_1,V_2,V_3;E)$ be a balanced
tripartite graph on $3M$ vertices so that each vertex is adjacent to
at least $\left(2/3-\epsilon'\right)M$ vertices in each of the
other classes.  If ${\cal T}_0$ is a partial cover by disjoint
triangles with $|{\cal T}_0|<M-\caa$, then we can find another partial
cover by disjoint triangles, $\cal T$ with $|{\cal T}|>|{\cal T}_0|$
and $|{\cal T}\setminus{\cal T}_0|\leq \cbb$, unless $G$ contains
three sets of size $M/3$ and pairwise have density less than
$\Delta_0$.
\label{lALMOSTCOVER}
\end{lemma}

\begin{proposition} If a reduced graph $G_r$ has two sets of size
$\ell/3$ and have density less than $\Delta_0$, then
some vertices can be added to the underlying graph induced by those
clusters so that it is two sets of size $\lfloor N/3\rfloor$ and have
density less than $\Delta$.
\label{pEXTREME}
\end{proposition}

Call these super-triangles $S(1), S(2), \ldots, S(\ell)$.  We put the
vertices in the remaining clusters into the appropriate leftover set.
Let the reduced graph involving the clusters of $S(1), S(2), \ldots,
S(\ell)$ be denoted $G_r$.  By Proposition~\ref{pSUPERREG}, at most
$2\epsilon_1 L'$ vertices can be removed from each cluster to obtain
$(\epsilon_3, \delta_3)$-super-regular pairs in the vertex-disjoint
triangular decomposition of $G_r$.  Furthermore,
Proposition~\ref{pSLICINGREG} guarantees that any edge in $G_r$ must
still correspond to an $\epsilon_3$-regular pair of density at least
$d_3$.
\begin{proposition} Given $\epsilon<1/4$, let $(S_i',S_j')$ for
$\{i,j\}\in {[3]\choose 2}$, be three $\epsilon$-regular pairs with
density at least $d$ and $|S_i'|=L'$ for $i=1,2,3$.  Some vertices can
be removed from each $S_i'$ to create $S_1$, $S_2$ and $S_3$ that form
three pairwise $(2\epsilon,d-3\epsilon)$-super-regular sets of size
$L\geq(1-2\epsilon)L'$.
\label{pSUPERREG}
\end{proposition}

\begin{proposition} Let $|X|=|Y|$, $X'\subseteq X$, $Y'\subseteq Y$,
$|X'|=|Y'|$ with $|X'|>\epsilon|X|$.  If $(X,Y)$ is
$\epsilon$-regular, then $(X',Y')$ is
$\max\left\{\left(\frac{|X|}{|X'|}\epsilon\right),
2\epsilon\right\}$-regular. 
\label{pSLICINGREG}
\end{proposition}

One cluster $y$ is {\bf reachable} from another, $x$, if there is a
chain of super-triangles, $T_1,\ldots,T_{2k}$ ($k\in\{1,2\}$) with $x$
an endpoint of $T_1$, and $y$ an endpoint of $T_{2k}$ with the added
condition that $T_{2i+1}$ and $T_{2i+2}$ ($i=0,\ldots,k-1$) share a
common edge and $T_{2i}$ and $T_{2i+1}$ ($i=1,\ldots,k-1$) share only
one common vertex.  Fix one super-regular super-triangle, $S(1)$.  The
set of all such triangles that connects some cluster to a cluster of
$S(1)$ is a {\bf structure}.  We would like to show that each cluster
in $G_r$ and $V_i$ is reachable from the cluster that is $S(1)\cap
V_i$.  If this is not possible, then Lemma~\ref{lREACHABLE} and
Proposition~\ref{pEXTREME} imply that $G$ must be in the extreme case.
\begin{lemma}[Reachability Lemma] 
In the reduced graph $G_r$, all clusters are reachable from other
clusters in the same class, unless some edges can be deleted from
$G_r$ so that the resulting graph obeys the minimum degree condition,
but is $\Delta_0$-approximately $\Theta_{3\times 3}(\ell/3)$.
\label{lREACHABLE}
\end{lemma}

So, suppose that every cluster is reachable from the appropriate
cluster of $S(1)$.  Consider some cluster $y$ and the structure that
connects it to $x$.  This structure contains clusters from at most 8
of the $S(i)$, not including $S(1)$ itself.  For any such structure,
$T_1,\ldots,T_{2k}$, find $\caa$ real triangles in each of the $T_i$,
for $i$ odd.  Note that if some $T$ is in more than one structure,
then there exist $\caa$ real triangles for each time that $T$ occurs
in a structure.  Do this for all possible structures, ensuring that
these real triangles are mutually disjoint and color these real
triangles red.  No cluster can possibly contain more than
$r=\threetimescaa\ell$ red vertices.  Thus, there are still $L-r$
uncolored vertices in each cluster, but
$L\geq[1-O(\epsilon_1)]\frac{N}{\ell}$, which goes to infinity as
$N\rto\infty$.  Proposition~\ref{pMAKETRIS} gives that finding these
red triangles is easy.
\begin{proposition}
Let $(X_1,X_2,X_3)$ be a triple with $|X_i|=L$ for $i=1,2,3$ and
each pair is $\epsilon$-regular with density $d>3\epsilon$.  Then,
there exist $(1-2\epsilon)L$ disjoint real triangles in the graph
induced by $(X_1,X_2,X_3)$. 
\label{pMAKETRIS}
\end{proposition}

This process of creating red triangles may result in an unequal number
of red vertices in the clusters of some of the $S(i)$'s.  Let $s_i$
denote the maximum number of red vertices in any one class of $S(i)$.
Pick a set of uncolored vertices of size $L-s_i$ in each class of
$S(i)$.  Proposition~\ref{pSLICINGREG} gives that the pairs of $S(i)$
are $(\epsilon',\delta')$-super-regular for some $\epsilon'$ and
$\delta'$.  Then, apply the Blow-up Lemma (Lemma~\ref{blowuplem}) to
get $L-s_i$ disjoint triangles among the uncolored vertices of $S(i)$.
Color these triangles blue.

Now, place the remaining uncovered vertices into the leftover sets.
Apply the Almost-covering Lemma (Lemma~\ref{lALMOSTCOVER}) to the
non-red vertices of $G$.  Each time this is applied, we may end up
destroying at most $\cbb$ of the blue triangles in order to create our
larger covering.  So, suppose that, at some point, there are less than
$(1-\delta_4)L+\caapcbb$ vertices remaining in some $S(i)$, then we
still apply the Almost-covering Lemma, but this time exclude vertices
in the blue triangles of $S(i)$ as well as red vertices.  There are at
most $\epsilon_5\ell$ of the $S(i)$'s that we may have to exclude in
this manner.

Color green any new triangles formed by using the Almost-covering
Lemma (Lemma~\ref{lALMOSTCOVER}).  There are at most $\threetimescaa$
uncolored vertices that remain after we are finished.  Let $x_1\in
V_1$ be an uncolored vertex.  We will show how to insert this vertex;
inserting the other vertices is similar.

The cluster containing $x_1$ has degree at least $2\delta_4L$ in at
least $(2/3-\alpha)\ell$ of the clusters in $V_2$ and $V_3$.
So, choose some $S(i)$ where $x_1$ is adjacent to at least $2\delta_4
L$ vertices in the $V_2$ and $V_3$ clusters of $S(i)$.  Color $x_1$
blue.  Now look at the structure that connects $S(i)$ to $S(1)$, and
call the triangles in this structure $T_1, \ldots, T_{2k}$.  Find a
triangle between the blue vertices of $T_{2k}$.  Color the edges and
vertices of this triangle red.  Next take one of the red triangles
from $T_{2k-1}$, uncolor its edges and color its vertices blue.
Continue in the same manner, adding a red triangle to $T_{2\xi}$ and
removing one from $T_{2\xi-1}$ for $\xi=k,k-1,\ldots,1$.  At the end
of this process, the same number of blue vertices are in each cluster
of each $S(j)$, except for one extra in $V(S(1))\cap V_1$.

Apply the same procedure to uncolored vertices in $V_2$ and $V_3$.
Now, the same number of blue vertices are in each $S(j)$, including
$S(1)$, which now has $\threetimescaa$ more blue vertices in each
class than before inserting the extra vertices.  Finally, apply the
Blow-up Lemma (the pairs are $(2\epsilon_3,\delta_4)$-super-regular)
to the blue vertices in each of the $S(j)$'s to create vertex-disjoint
blue triangles that involve all of the blue vertices.  So, the red,
green and blue triangles are vertex-disjoint and cover all vertices of
$G$. $\EOP$

\subsection{Proofs of Propositions}
\noindent{\bf Proof of Proposition~\ref{pEXTREME}.}  This is
immediate from the fact that the density of any pair of clusters
nonadjacent in $G_r$ is at most $d_1+2\epsilon_1$ and from the fact
that $\Delta_0\ll\Delta$. $\EOP$



\noindent {\bf Proof of Proposition~\ref{pSUPERREG}.} Let $T$ be the
subset of $S_1'$ consisting of vertices with degree at most
$(d-\epsilon)L'$ in $S_2'$. Clearly $d(T, S_2')\leq d-\epsilon$.
But, if $|T|>\epsilon L'$, then $d(T,S_2')>d-\epsilon$, a
contradiction.  So, $|T|\leq\epsilon L'$.  We then have at least
$(1-2\epsilon)L'$ vertices in $S_1'$ that have degree at least
$(d-\epsilon)L'$ in both $S_2'$ and $S_3'$.  Call that set $S_1$ and
similarly define $S_2$ and $S_3$.  Proposition~\ref{pSLICINGREG}
gives that these sets are pairwise $2\epsilon$-regular if
$\epsilon<1/4$, then the proposition is proven. $\EOP$

\noindent A proof of Proposition~\ref{pSLICINGREG} is
straightforward and left to the reader. 

\noindent{\bf Proof of Proposition~\ref{pMAKETRIS}.} We apply
Proposition~\ref{pSUPERREG} to the triple $(X_1,X_2,X_3)$
to get a triple $(X_1',X_2',X_3')$ such that each pair is
$(2\epsilon,d-3\epsilon)$ super-regular each on
$L^{*}\geq(1-2\epsilon)L$ vertices.  We then apply the Blow-up Lemma
(Lemma~\ref{blowuplem}) to $(X_1',X_2',X_3')$ getting our $L^{*}$
vertex-disjoint triangles. $\EOP$.

\subsection{Proof of the Almost-covering Lemma
(Lemma~\ref{lALMOSTCOVER})}

Given the constants $\epsilon'\ll\Delta'\ll\Delta_0$, let ${\cal T}_0$
be as in the statement of the lemma.  Denote $U_1$, $U_2$ and $U_3$ as
the portions of $V_1$, $V_2$, and $V_3$, respectively, left uncovered
by ${\cal T}_0$. Let $U=U_1+U_2+U_3$.  We want to show that if $|U|>9$
then the covering can be expanded unless $G$ contains three sets of
size $M/3$ which pairwise have density less than $\Delta_0$.  (We
always assume that $M$ is divisible by $3$.)

Thus, assume that $U$ contains at least four vertices in each class.
We want to show that there are at least three disjoint edges, one
between each class.  Let $x_1\in U_1$ and $x_2\in U_2$ with
$x_1\not\sim x_2$ then it will be possible to exchange these vertices
with the vertices of $\cal T$ that maintains or increases the number
of disjoint triangles, uses no other vertices in $U$ and places an
edge between $U_1$ and $U_2$.

By assumption, both
$\left|N_{V_2\setminus U_2}(x_1)\right|\geq
\left(2/3-\epsilon'\right)M$ and $\left|N_{V_3\setminus U_3}
(x_1)\right|\geq (2/3-\epsilon')M$.  This implies that there are at
least $(1/3-2\epsilon')M$ triangles, $T$, in ${\cal T}_0$ so that
$x_1$ is adjacent to both the $V_2$ and $V_3$ vertices in $T$.

Let
\begin{eqnarray*}
   A_1 & := & \left\{x\in V_1 : T\in{\cal T}_0, V(T)=\{x,y,z\}, 
   x_1\sim y, \mbox{ and } x_1\sim z\right\} \\
   A_2 & := & \left\{y\in V_2 : T\in{\cal T}_0, V(T)=\{x,y,z\}, 
   x_1\sim y, \mbox{ and } x_1\sim z\right\} \\
   A_3 & := & \left\{z\in V_3 : T\in{\cal T}_0, V(T)=\{x,y,z\}, 
   x_1\sim y, \mbox{ and } x_1\sim z\right\} 
\end{eqnarray*}
Simply, $A_1$ is the set of all vertices so that $x_1$ can be
exchanged with such a vertex so as to leave the number of triangles in
${\cal T}_0$ unchanged.  The sets $A_2$ and $A_3$ are the vertices in
the other classes that correspond to the triangles in ${\cal T}_0$
with vertices in $A_1$.  Clearly, $|A_1|=|A_2|=|A_3| \geq
\left(1/3-2\epsilon'\right)M$.

Consider $x_2$.  Define $B_1$, $B_2$ and $B_3$ in a similar manner
so that $x_2$ can be exchanged with each of the vertices of $B_2$.  We
will show that the intersection of $A_1$ and $B_1$ is empty.  If there
is a triangle $\{x,y,z\}\in{\cal T}_0$ such that $x\in A_1\cap B_1$,
then $x$ and $x_1$ can be exchanged in order to obtain a covering with
the same number of triangles but with an edge in $(U_1,U_2)$.

The pair $(A_1,B_2)$ is void of edges.  If not, then both $x_1$ and
$x_2$ can be exchanged with the endvertices of that edge.  The number
of triangles does not change, but there will be an edge between $U_1$
and $U_2$.  Now let $C_i=V_i({\cal T}_0)\setminus(A_i\cup B_i)$, for
$i=1,2,3$.  Clearly, $|C_1|=|C_2|=|C_3|\leq (1/3+4\epsilon')M$.  But,
since no vertex in $A_1$ can be adjacent to a vertex in $U_2$ and
$(A_1,B_2)$ is void, then $|C_2|\geq (1/3-2\epsilon')M$.

\begin{figure}
\centerline{\epsfbox{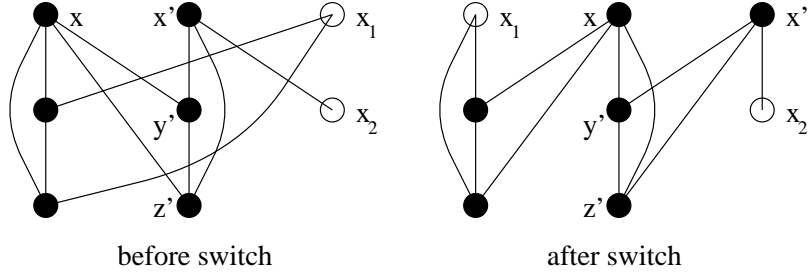}}
\caption{$(A_1,C_3)$ not sparse}
\label{figure2}
\end{figure}

Since $(A_1,B_2)$ is void, each vertex in $A_1$ must be adjacent to at
least $(1/3-3\epsilon')M$ vertices in $C_2$.  Therefore, if there
exists some vertex $x\in A_1$ adjacent to more than $7\epsilon'M$
vertices in $C_3$, then there is a triangle $\{x',y',z'\}$ such that
$x\sim y',z'$.  Thus, according to Figure~\ref{figure2}, $x$, $x'$ and
$x_1$ could be moved so that there exists a ${\cal T}_0$ of the same
size with an edge in $(U_1,U_2)$.

Similarly, $(B_2,C_3)$ must be sparse.  Therefore, the triple
$(A_1,B_2,C_3)$ has sets of size $4\epsilon'$-approximately $M/3$ with
pairwise density at most $7\epsilon'M/|C_3|<\Delta'$.  The same
procedure can be applied to $(U_1,U_3)$ and then $(U_2,U_3)$ to create
6 edges, $e_1,e_2\in (U_1,U_2)$, $f_1,f_3\in (U_1,U_3)$ and
$g_2,g_3\in (U_2,U_3)$ that are disjoint.  Let this new partial
triangular cover be ${\cal T}_2$. Note that $|{\cal T}_0\setminus
{\cal T}_2|\leq 12$ but $|{\cal T}_2|=|{\cal T}_0|$.

Given the edges $e_1$, $e_2$, $f_1$, $f_3$, $g_2$, $g_3$ in $U$,
redefine the ``$A$'', ``$B$'', and ``$C$'' sets.  Let
\begin{eqnarray*}
   A_i & = & \left\{x_i\in V_i : \{x_1,x_2,x_3\}\in{\cal T}_2
   \mbox{ and }x_1\sim g_i\right\}\qquad i=2,3 \\
   B_i & = & \left\{x_i\in V_i : \{x_1,x_2,x_3\}\in{\cal T}_2
   \mbox{ and }x_2\sim f_i\right\}\qquad i=1,3 \\
   C_i & = & \left\{x_i\in V_i : \{x_1,x_2,x_3\}\in{\cal T}_2
   \mbox{ and }x_3\sim e_i\right\}\qquad i=1,2
\end{eqnarray*}
The sets $|A_i|,|B_i|,|C_i|\geq (1/3-2\epsilon')M$ for all
relevant $i$.  This is the case because the neighborhood of each edge
is of this size and these neighborhoods must be entirely within
$V({\cal T}_2)$, otherwise ${\cal T}_2$ is not maximal.

We wish to show that these sets are disjoint.  Suppose, without loss
of generality, $x_1\in B_1\cap C_1$ so that $\{x_1,x_2,x_3\}\in {\cal
T}_2$.  Then $x_2$ and $f_1$ form a triangle and $x_3$ and $e_1$ form
a triangle -- giving that there exists a $\cal T$ of size larger than
${\cal T}_2$.
As a result, $|B_1\cup C_1|,|A_2\cup C_2|,|A_3\cup B_3|\geq
(2/3-4\epsilon')M$.

\begin{figure}
\centerline{\epsfbox{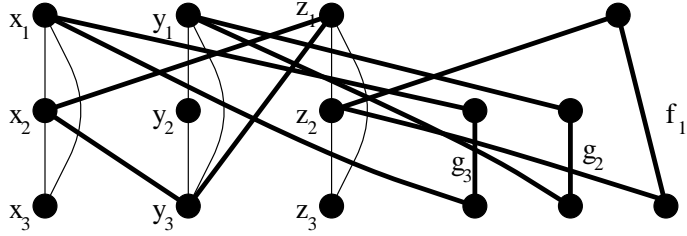}}
\caption{A triangle in $(B_1,A_2,A_3)$}
\label{figure4}
\end{figure}

Further, there can be no triangle in the triple $(B_1\cup C_1, A_2\cup
C_2, A_3\cup B_3)$.  We will just show one example; suppose there is a
triangle $T$ in $(B_1,A_2,A_3)$.  Then there are
$\{x_1,x_2,x_3\},\{y_1,y_2,y_3\},\{z_1,z_2,z_3\}\in{\cal T}_2$ such
that $x_1\sim g_2$, $y_1\sim g_3$, $z_2\sim f_1$ and
$T=\{z_1,x_2,y_3\}$.  If $x_1=y_1$, then $\{x_1,x_2,x_3\}$ and
$\{z_1,z_2,z_3\}$ can be replaced with the triangle formed by $x_1$
and $g_2$, the triangle formed by $z_2$ and $f_1$ and $T$ itself.  If
$x_1\neq y_1$, then we can replace the $x$, $y$ and $z$ triangles with
the triangles formed by $x_1$ and $g_2$, by $y_1$ and $g_3$ and by
$z_2$ and $f_1$ as well as $T$.  See Figure~\ref{figure4}.
Thus, if there were a triangle in $(B_1\cup C_1, A_2\cup C_2, A_3\cup
B_3)$, then a $\cal T$ could be found such that $|{\cal T}|>|{\cal
T}_0|$ and $|{\cal T}_0\setminus{\cal T}|\leq\cbb$.  If there is no
triangle in $(B_1\cup C_1, A_2\cup C_2, A_3\cup B_3)$, then
Proposition~\ref{pTHETA} gives that this subgraph is
$\Delta_0$-approximately $\Theta_{3\times 2}(M/3)$,
which means $G$ contains 3 sets of size $M/3$ with pairwise density at
most $\Delta_0$. $\EOP$

\subsection{Proof of the Reachability Lemma
(Lemma~\ref{lREACHABLE})} 
Let us be given constants
\[ \epsilon_2\ll\Delta'\ll\Delta''\ll\Delta'''\ll\Delta_0 . \]

In order to prove the lemma, we distinguish two triangles, call them
$S(1)=\{x_1(1), x_2(1),x_3(1),\}$ and
$S(\ell)=\{x_1(\ell),x_2(\ell),x_3(\ell)\}$ and suppose $x_1(\ell)$ is
not reachable from $x_1(1)$.  We will show that edges can be deleted
from $G_r$ so that the minimum degree condition holds and the
resulting graph is $\Delta_0$-approximately $\Theta_{3\times
3}(\ell/3)$.  Every cluster is adjacent to at least
$(2/3-\epsilon_2)\ell$ clusters in each of the other classes.  Let
\begin{eqnarray*}
   A_{i,1} & := & \left[N(x_1(\ell))\setminus N(x_1(1))\right]\cap 
   V_i \\
   A_{i,3} & := & \left[N(x_1(1))\setminus N(x_1(\ell))\right]\cap 
   V_i \qquad\qquad i=2,3 \\
   A_{i,2} & := & V_i\setminus\left(A_{i,1}\cup A_{i,3}\right)
\end{eqnarray*}
Observe that $(1/3-2\epsilon_2)\ell\leq |A_{i,1}|,|A_{i,3}|\leq
(1/3+\epsilon_2)\ell$ for $i=2,3$.

If there is an edge in $(A_{2,2},A_{3,2})$, then $x_1(\ell)$ must be
reachable from $x_1(1)$.  Thus, it must be that
$d(A_{2,2},A_{3,2})=0$.  Combining the information, it must be true
that $|A_{i,j}|\in\left((1/3-4\epsilon_2)\ell,
(1/3+4\epsilon_2)\ell\right)$ for $i\in\{2,3\}$ and $j\in\{1,2,3\}$.
Define the sets $A_{1,1}$ and $A_{1,3}$ by first letting
\[ A_{1,1}\cup A_{1,3} :=\left\{v\in V_1 : \exists
   i\in\{2,3\}\mbox{ s.t. }\deg_{A_{i,2}}(v)\geq 2\Delta'\ell\right\}
   \; . \] 

Suppose $v\in A_{1,1}\cup A_{1,3}$ with $\deg_{A_{i,2}}(v)\geq
2\Delta'\ell$ and $\deg_{A_{i',1}}(v),
\deg_{A_{i',3}}(v)\geq\Delta''\ell$, where
$\{i'\}=\{2,3\}\setminus\{i\}$.
Then there exists an edge in $(A_{i,2},A_{i',3})$ that is adjacent to
both $x_1(1)$ and $v$.  Also, there exists another edge in
$(A_{i,2},A_{i',1})$ that is adjacent to both $v$ and $x_1(\ell)$.
This makes $x_1(\ell)$ reachable from $x_1(1)$ by a chain of 4
triangles.

Suppose $v\in A_{1,1}\cup A_{1,3}$ and $v$ is adjacent to less than
$\Delta''\ell$ vertices in $A_{i',1}$ but is adjacent to more than
$\Delta''\ell$ vertices in $A_{i,1}$.  In this case, there exists an
edge in $(A_{i,2},A_{i',3})$ that is adjacent to both $x_1(1)$ and
$v$.  There also exists an edge in $(A_{i,1},A_{i',2})$ that is
adjacent to both $v$ and $x_1(\ell)$.  With the above suppositions
about the degree of $v$, $x_1(\ell)$ is reachable from $x_1(1)$ by a
chain of 4 triangles.  Therefore, each vertex either is adjacent to
less than $\Delta''\ell$ vertices in both $A_{2,4}$ and $A_{3,4}$
(call these vertices $A_{1,1}$) or is adjacent to less than
$\Delta''\ell$ vertices in both $A_{2,1}$ and $A_{3,1}$ (call these
vertices $A_{1,3}$).  This gives that $d(A_{1,1},A_{i,3})<\Delta'''$
and $d(A_{1,3},A_{i,1})<\Delta'''$ for $i=2,3$.

Because of the minimum degree condition,
$|A_{1,1}|,|A_{1,3}|<(1/3+2\epsilon_2)\ell$. 
Define $A_{1,2}$ to be those vertices adjacent to less than
$2\Delta'\ell$ vertices in both $A_{2,2}$ and $A_{3,2}$.  It must be
true that $V_1=A_{1,1}\cup A_{1,2}\cup A_{1,3}$ with all sets being
disjoint, because the definition of $A_{1,1}\cup A_{1,3}$ gives that
all vertices not in those sets must be in $A_{1,2}$.  From before,
$|A_{1,2}|<(1/3+2\epsilon_2)\ell$ and $d(A_{1,2},A_{i,2})<\Delta'''$
for $i=2,3$.  Summarizing, $|A_{i,j}|\in\left((1/3-O(\epsilon_2))\ell,
(1/3+O(\epsilon_2))\ell\right)$ for all $i$ and $j$.  Furthermore,
$d(A_{1,1},A_{i,3})<\Delta'''$ and $d(A_{1,3},A_{i,1})<\Delta'''$ for
$i=2,3$ and, as we just showed, $d(A_{1,2},A_{i,2})<\Delta'''$ for
$i=2,3$.

What remains is to show that one of the pairs $(A_{2,1},A_{3,1})$ or
$(A_{2,3},A_{3,3})$ is pairwise sparse.  Note that if one is pairwise
sparse, we might as well allow the other to be pairwise sparse, since
extra edges only help.  If it is not the case, then
$d(A_{2,1},A_{3,1})\geq\Delta'''$ and
$d(A_{2,3},A_{3,3})\geq\Delta'''$ are both not sparse.  There exists
an edge $e$ in $(A_{2,1},A_{3,1})$ that is adjacent to many vertices
in $A_{1,2}$ as well as an edge $f$ in $(A_{2,3},A_{3,3})$ that is
adjacent to many vertices in $A_{1,2}$.  Thus, there exists a vertex
$v$ that is adjacent to both edges.  Since $f$ is adjacent to both
$x_1(1)$ and $v$ and $e$ is adjacent to both $v$ and $x_1(\ell)$ we
have that $x_1(\ell)$ is reachable from $x_1(1)$.  Therefore, if
$x_1(\ell)$ is not reachable from $x_1(1)$, we must have that
$d(A_{2,1},A_{3,1}), d(A_{2,3},A_{3,3})\ll\Delta_0$. $\EOP$

\section{The Extreme Tripartite Theorem}
\label{sTRIEXTREME}
\subsection{Statement of the Theorem}

Theorem~\ref{tTRIFUZZY} leaves the extreme case, which we
consider in Theorem~\ref{tTRIEXTREME}.

\begin{theorem} Given $\Delta\ll 1$, let $G=(V_1,V_2,V_3;E)$ be
a balanced tripartite graph on $3N$ vertices such that each vertex is
adjacent to at least $(2/3)N$ vertices in each of the other classes.
Furthermore, let $G$ have three sets with size $N/3$ and pairwise
density at most $\Delta$.  Then, if $N$ is large enough, either $G$
can be covered with triangles or $G$ is $\Gamma_3(N/3)$.
\label{tTRIEXTREME}
\end{theorem}

\subsection{Proof of the Theorem}
Assume that $G$ is minimal.  That is, no edge of $G$ can be deleted so
that the minimum degree condition still holds.  We will prove that for
minimal $G$, either $G$ can be covered with triangles or
$G=\Gamma_3(N/3)$.  With that proven, it suffices to show that adding
any edge to $\Gamma_3(N/3)$ will allow the resultant graph to be
covered with triangles -- this will be discussed in
Section~\ref{sGAMMACASE}.  Begin with the usual sequence of constants:
\[ \Delta\ll\Delta_1\ll\Delta_2\ll\eta\ll\theta-\frac{3}{4} \] 
for some $\theta$, $3/4<\theta<1$.  Let $t:=N/3$ with $N$ divisible by
3.

Let the sets of size $t$ mentioned in the theorem be designated $A_i$,
with $A_i\subset V_i$ for $i=1,2,3$.  Let $B_i:=V_i\setminus A_i$ for
$i=1,2,3$.  For each $i\in\{1,2,3\}$, let $A_i'$ be the vertices that
are adjacent to at least $(1+\theta)t$ vertices in $B_j$ for each
$j\neq i$.  Let $B_i'$ be the vertices that are adjacent to at least
$(1/2)(1+\theta) t$ vertices in $A_j$ for each $j\neq i$.
Furthermore, let $C_i'=V_i\setminus\left(A_i'\cup B_i'\right)$.  The
key feature of each $c\in C_i'$ is that there is a $j\neq i$ such that
$c$ is adjacent to at least $(1-\theta)t$ vertices in $A_j$.  Let us
compute $|A_i'|$ and $|B_i'|$ for $i=1,2,3$.  Proposition~\ref{pSIZES}
restricts the sizes of these sets.
\begin{proposition}
If $\Delta\ll\Delta_1$, then for all $i\in\{1,2,3\}$,
\begin{eqnarray*}
   |A_i'| & \in & \left((1-\Delta_1)t, (1+\Delta_1)t\right) \\
   |B_i'| & \in & \left((2-\Delta_1)t, (2+\Delta_1)t\right)
\end{eqnarray*}
Furthermore, for $i=1,2,3$, $|A_i\setminus A_i'|,|B_i\setminus
B_i'|\leq\Delta_1 t$.
\label{pSIZES}
\end{proposition}

The key lemma for this proof is Lemma~\ref{lKEY}.
\begin{lemma}
Let $\Delta_1\ll\Delta_1'\ll\Delta_1''\ll\Delta_1'''\ll\Delta_2'\ll
\Delta_2\ll\theta-\frac{3}{4}$ for some constant $\theta>3/4$.  Let
$G=(V_1,V_2,V_3;E)$ be a balanced tripartite graph on $9t'$ vertices
with each vertex adjacent to at least $(2-\Delta_1')t'$ vertices in
each of the other vertex classes.  Suppose further that we have sets
$A_i'$ of size $\Delta_1''$-approximately $t'$ such that for all $a\in
A_i'$, $\deg_{V_j\setminus A_j'}(a)\geq (1+\theta)t'$ for all $j\neq
i$.  Furthermore, let $d(A_i',A_j')<\Delta_1'''$, $\forall
\{i,j\}\in{[3]\choose 2}$ and let each $v\in V_i\setminus A_i'$ have
the property that there is a $j\neq i$ such that $\deg_{A_j'}(v)\geq
(1-\theta-\Delta_1'')t'$.  If $G$ is minimal and cannot be covered
with triangles then either
\begin{enumerate}
   \item $|A_1'|+|A_2'|+|A_3'|>3t$, \label{KEY1} 
   \item $G$ is $\Delta_2'$-approximately $\Gamma_3(t')$, or
\label{KEY2} 
   \item $G$ is $\Delta_2'$-approximately $\Theta_{3\times
3}(t')$. \label{KEY3}
\end{enumerate}
\label{lKEY}
\end{lemma}

We make adjustments according to whether or not
$|A_1'|+|A_2'|+|A_3'|\leq 3t$.  It is true that, $3(1-\Delta_1)t\leq
|A_1'|+|A_2'|+|A_3'|\leq 3(1+\Delta_1)t$.  If
$|A_1'|+|A_2'|+|A_3'|\leq 3t$, then apply Lemma~\ref{lKEY} to $G$.
Thus, $G$ can be covered with triangles unless $G$ is
$\Delta_2$-approximately $\Gamma_3(t)$ or $G$ is
$\Delta_2$-approximately $\Theta_{3\times 3}(t)$.

If $|A_1'|+|A_2'|+|A_3'|>3t$, then we want to create a matching of
size $|A_1'|+|A_2'|+|A_3'|-3t$ in $(A_1',A_2',A_3')$.  After finding
the matching, find common neighbors in $(B_1'\cup C_1',B_2'\cup
C_2',B_3'\cup C_3')$ and remove those disjoint triangles so that
Lemma~\ref{lKEY} can be applied to the remaining graph.  Each vertex
in $A_i'$ is adjacent to at least $\max\{|A_j'|-t,0\}$ vertices in
$A_j'$, for all distinct $i$ and $j$.  Thus, we can create matchings
sequentially in $(A_i',A_j')$, for all pairs $(i,j)$ so that they do
not coincide and together they exclude exactly $t$ vertices in each of
$A_1'$, $A_2'$ and $A_3'$ that are larger than $t$.  The details are
left to the reader.



Thus, $G$ can be covered with triangles unless $G$ is either
$\Delta_2$-approximately $\Gamma_3(t)$ (Section~\ref{sGAMMACASE}) or
$\Delta_2$-approximately $\Theta_{3\times 3}(t)$
(Section~\ref{sTHETACASE}).

\subsection{$G$ is $\Delta_2$-approximately $\Gamma_3(t)$ \\ (Case
(\ref{KEY2}) of Lemma~\ref{lKEY})}
\label{sGAMMACASE}
\begin{figure}
\centerline{\epsfbox{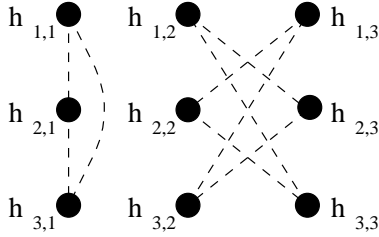}}
\caption{Diagram of $\Gamma_3$.  The dotted lines correspond to non-edges.}
\label{efigureGAMMA}
\end{figure}

Let the sets $A_{i,j}$, $i,j=1,2,3$ be as the $h_{i,j}$ in
Figure~\ref{efigureGAMMA}.  
Note that the figure depicts the non-edges of this graph.  Each row of
vertices corresponds to a vertex class and the dotted lines correspond
to non-edges.  Given $\Delta_2\ll\Delta_3\ll\Delta_4\ll\eta$, our goal
is to modify the sets $A_{i,j}$ to form sets $\tilde{A}_{i,j}$.  The
triangles will come from each of the following:
\begin{eqnarray*}
   (\tilde{A}_{1,1},\tilde{A}_{2,2},\tilde{A}_{3,2}) &
   (\tilde{A}_{2,1},\tilde{A}_{1,2},\tilde{A}_{3,2}) &
   (\tilde{A}_{3,1},\tilde{A}_{1,2},\tilde{A}_{2,2}) \\
   (\tilde{A}_{1,1},\tilde{A}_{2,3},\tilde{A}_{3,3}) & 
   (\tilde{A}_{2,1},\tilde{A}_{1,3},\tilde{A}_{3,3}) & 
   (\tilde{A}_{3,1},\tilde{A}_{1,3},\tilde{A}_{2,3}) .
\end{eqnarray*}
The triangles will receive one of 6 labels $(i;j)$, for
$i\in\{1,2,3\}$ and $j\in\{2,3\}$.  A triangle with the label $(i;j)$
will be in the triple $(\tilde{A}_{i,1}, \tilde{A}_{i_2,j},
\tilde{A}_{i_3,j})$, where $i_2,i_3$ are distinct indices in
$\{1,2,3\}\setminus\{i\}$.

Define $A_{i,j}'$ to be the set of ``typical'' vertices in $A_{i,j}$.
That is, if $\{h_{i_1,j_1},h_{i_2,j_2}\}$ is a non-edge in $\Gamma_3$,
then each vertex in $A_{i_1,j_1}'$ is adjacent to less than $\eta t$
vertices in $A_{i_2,j_2}$.  Let $C_i=V_i\setminus\left(A_{i,1}'\cup
A_{i,2}'\cup A_{i,3}'\right)$, for $i=1,2,3$.  Since
$\Delta_2\ll\Delta_3$, $|A_{i,j}\setminus A_{i,j}'|<\Delta_3 t$.  We
will make the sets $A_{i,1}'$ into sets $A_{i,1}''$ of size $t$, for
$i=1,2,3$.  If there is some $|A_{i,1}'|>t$, then we find a matching
in $(A_{1,1}',A_{2,1}',A_{3,1}')$ of size
$\sum_{i=1}^3\max\{|A_{i,1}'|-t,0\}$ similar to the one we constructed
above.  Color this matching red and for $|A_{i,1}'|>t$, take
$|A_{i,1}'|-t$ red edges and remove the $V_i$ endvertices that are in
$A_{i,1}'$ and add them to one of $A_{i,2}'$ or $A_{i,3}'$, whichever
has size smaller than $t$.  This creates sets $A_{i,1}''$ of size at
most $t$, for $i=1,2,3$.  The endvertices of this red edge will
receive label $(i;j)$ if one of its vertices is added to $A_{i,j}'$.

Suppose that $|A_{i,1}'|<t$.  Then find vertices in either $A_{i,2}'$
or $A_{i,3}'$, color them green and add them to $A_{i,1}'$ to form
$A_{i,1}''$.  To show that these green vertices will act as $A_{i,1}'$
vertices, suppose, without loss of generality, $v$ is a green vertex
added to $A_{1,1}'$.  Observe that $v$ must be adjacent to at least
$(1-\eta)t$ vertices in either $A_{2,2}$ and $A_{2,3}$ or $A_{2,3}$
and $A_{3,3}$.  Thus, if we move a green vertex from $A_{i,j}'$, it
will receive the label $(i;j)$.  The resulting sets $A_{i,1}''$ are of
size exactly $t$, so let them be renamed $\tilde{A}_{i,1}$, $i=1,2,3$.

Now we want to show that vertices in $C_i'$ behave like vertices in
either $A_{i,2}'$ or $A_{i,3}'$.  Let $c\in C_1'$ and, without loss of
generality, show that $c$ can be added to $A_{1,2}'$.  There exists an
$i\in\{2,3\}$ such that $c$ is adjacent to at least $\eta t$ vertices
in $A_{i,1}$.  If $c$ is adjacent to at least $\eta t$ vertices in
$A_{5-i,2}$, then $c$ can receive the label $(i;2)$.  Otherwise $c$
can receive the label $(5-i;2)$.  Color the $C_i'$ vertices green and
add them to either $A_{i,2}'$ or $A_{i,3}'$ (the smaller of the two)
to form $A_{i,2}''$ and $A_{i,3}''$.

Unfortunately, one of the sets $A_{i,2}''$ or $A_{i,3}''$ might be of
size more than $t$.  In order to create sets of size $t$, let us
suppose without loss of generality, that both $|A_{1,2}''|>t$ and
$A_{1,2}''$ is the largest from among $A_{1,2}''$, $A_{1,3}''$,
$A_{2,2}''$ and $A_{2,3}''$.  Let $\tau=|A_{1,2}''|-t$ and observe
that $A_{1,2}''=A_{1,2}'$.  Let $q=|A_{1,2}''|-|A_{2,3}''|\leq 2\tau$
because $|A_{1,2}''|\geq|A_{2,2}''|$.  Let $W\subset A_{2,3}''\cap
A_{2,3}'$,
\[ \left(|A_{1,2}''|-t\right)|W|\leq e(W,A_{1,2}'')\leq
   (\gamma+2\Delta_2) t\left|N_{A_{1,2}''}(W)\right| . \] 
So, $|N_{A_{1,2}''}(W)|\gg 2\tau$,
provided $|W|$ is not too small, and there exists a matching of size
$q$ in $(A_{1,2}'',A_{2,3}''\cap A_{2,3}')$.  Color this matching
blue.

If $|A_{2,2}''|>t$, take $|A_{1,2}''|-|A_{2,2}''|$ blue vertices from
$A_{1,2}''$ and add them to $A_{1,3}''$.  Also, take $|A_{2,2}''|-t$
edges in $(A_{1,2}'',A_{2,2}'')$, color them blue and add their
vertices to $A_{1,3}''$ and $A_{2,3}''$.  Such blue edges will be in
triangles with label $(3;3)$.  If $|A_{2,2}''|\leq t$, then take
$t-|A_{2,2}''|$ blue vertices from $A_{2,3}''$ and add them to
$A_{2,2}''$.  The endvertices of these blue edges will be in triangles
labeled $(3;2)$.  For the remaining blue edges, take
$\left|A_{1,2}''\right|-t$ of the vertices from $A_{1,2}''$ and add
them to $A_{1,3}''$.  The endvertices of these blue edges will be in
triangles labeled $(3;3)$.

It may be necessary to find a similar matching in $(V_2,V_3)$ if either
$|A_{3,2}''|>t$ or $|A_{3,3}''|>t$.  It is easy to see that we can do
so without using any of the other colored vertices by choosing a $W$
that excludes blue vertices.  The sets that result from moving the
vertices of blue edges are of size exactly $t$, so denote them
$\tilde{A}_{i,j}$, for $i=1,2,3$ and $j=2,3$.

Recall that $\Gamma_3(t)$ cannot be covered with triangles if $t$ is
odd.  A similar dilemma must also be resolved in this case.  Suppose
$t$ is odd.  Our goal is to find three triangles in $G$ such that each
vertex is from a different $\tilde{A}_{i,j}$.  Call these {\bf parity
triangles}.  To find them, we look for an edge in
$(\tilde{A}_{1,2}\cup\tilde{A}_{2,2}\cup\tilde{A}_{3,2},
\tilde{A}_{1,3}\cup\tilde{A}_{2,3}\cup\tilde{A}_{3,3})$ with a common
neighbor in $\tilde{A}_{1,1}\cup\tilde{A}_{2,1}\cup\tilde{A}_{3,1}$.

If there is such an edge among uncolored vertices, then there are many
neighbors in $\tilde{A}_{1,1}\cup\tilde{A}_{2,1}\cup\tilde{A}_{3,1}$.
If any vertex was colored, then by re-examining the process by
which it was constructed we see that it is possible to create such a
triangle.  For example, if there is a red edge in
$(\tilde{A}_{2,1},\tilde{A}_{1,2})$, then we can find a common
neighbor in $\tilde{A}_{3,3}$.  In any case, if the desired triangle
is found, remove it, along with two other triangles so that the
resulting ``$\tilde{A}$'' sets are of size $t-1$.  If a parity
triangle cannot be created in this way, then $G$ contains no colored
vertices.  In that case, if there is an edge in the graph that is
induced by $\tilde{A}_{1,1}\cup\tilde{A}_{2,1}\cup\tilde{A}_{3,1}$,
then it is easy to create the parity triangles.  Otherwise,
$G=\Gamma_3(t)$ and the theorem would be proven.

Therefore, suppose that the remaining $\tilde{A}_{i,j}$ sets are of
the same even cardinality.  Partition each $\tilde{A}_{i,j}$ uniformly
at random into two equally-sized pieces.  Each piece will receive one
of six labels $(i;j)$, for $i=1,2,3$ and $j=2,3$.  For
$i\in\{1,2,3\}$, $\tilde{A}_{i,1}$ will be partitioned into one set
labeled $(i;2)$ and the other labeled $(i;3)$.  For $i\in\{1,2,3\}$
and $j\in\{2,3\}$, $\tilde{A}_{i,j}$ will be partitioned into one set
labeled $(i_1;j)$ and the other labeled $(i_2;j)$, where $i,i_2,i_3$
are distinct members of $\{1,2,3\}$.

The triangle cover will only consist of triangles with vertices in
pieces with the same label.  Each of the colored vertices corresponds
to at least one of the two labels, but not necessarily both.  For
example, if there is a red edge in
$(\tilde{A}_{2,1},\tilde{A}_{1,2})$, then we want to ensure that each
of its endvertices are in pieces labeled $(2;2)$.  So, it may be
necessary to exchange colored vertices in one piece with uncolored
ones in the other piece.  A total of at most $\Delta_4 t$ vertices
will be so exchanged in any $\tilde{A}_{i,j}$.  

The covering by triangles can be completed by taking each piece that
has the same label and covering the corresponding triple with
triangles.  Consider, for example, the vertices in $\tilde{A}_{1,1}$,
$\tilde{A}_{2,2}$ and $\tilde{A}_{3,2}$ that carry the label $(1;2)$.
For simplicity, call them $S_1$, $S_2$ and $S_3$, respectively.

Any green vertex $v\in S_1$ is adjacent to at least $\eta t$ vertices
in both $A_{2,2}$ and $A_{3,2}$.  Thus, it is adjacent to at least
$(\eta-2\Delta_3)t$ vertices in both $\tilde{A}_{2,2}$ and
$\tilde{A}_{3,2}$.  Since the $S_i$ were chosen at random, Stirling's
inequality (see Corollary~\ref{cSTIRLING} in Section~\ref{sSTIRLING})
gives that $v$ is adjacent to at least
$(\eta-O(\Delta_4))(t/2)$ uncolored and unexchanged vertices in both
$S_2$ and $S_3$.  Since all but $O(\Delta_2)t$ of the vertices in
$A_{2,2}$ have degree at least $(1-2\Delta_2)t$ in $A_{3,3}$,
Stirling's inequality again gives that there exists an edge among the
uncolored and unexchanged vertices of $(N(v)\cap S_2,N(v)\cap S_3)$.
Do this for all the green vertices in order to get disjoint green
triangles.

The red and blue edges are even easier.  For example, since each
endvertex of a red edge in $(S_1,S_2)$ is adjacent to at least
$(1-\eta)t$ vertices in $A_{3,2}$, we can find a common vertex among
the uncolored and unexchanged vertices of $S_3$.  So, extend the
colored edges to find red and blue triangles disjoint from each other
and from the green triangles.  Finally each uncolored vertex in $S_i$
that was ``exchanged'' has degree at least
$(1-\eta-O(\Delta_3))(t/2)$ in each of the $S_j$, $\forall j\neq i$.
Put these in black triangles disjoint from each other and from other
colored triangles.  Let there be $t'$
$\geq(1-O(\Delta_4))(t/2)$ uncolored vertices remaining in each class.
Call them $S_i'\subset S_i$, for $i=1,2,3$.  Since $\Delta_4\ll\eta$,
each vertex in $S_i'$ is adjacent to at least $(3/4)t'$ vertices in
each of the $S_j$, $\forall j\neq i$.  Proposition~\ref{pEASYCASE}
finishes the covering and the proof of this case.

\subsection{$G$ is $\Delta_2$-approximately $\Theta_{3\times 3}(N/3)$ \\
(Case (\ref{KEY3}) of Lemma~\ref{lKEY})}
\label{sTHETACASE}
\begin{figure}
\centerline{\epsfbox{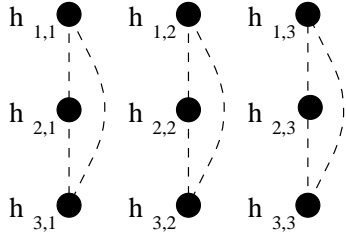}}
\caption{Diagram of $\Theta_{3\times 3}$.  The dotted lines correspond
to non-edges.} 
\label{efigureTHETA}
\end{figure}

Let the sets $A_{i,j}$, $i,j=1,2,3$ be as the $h_{i,j}$ in
Figure~\ref{efigureTHETA}.  Note that the figure depicts the non-edges
of this graph.  Each row of vertices corresponds to a vertex class and
the dotted lines correspond to non-edges.  
Our goal is again to modify the sets $A_{i,j}$ to form sets
$\tilde{A}_{i,j}$.  Our triangles will come from
$\left(\tilde{A}_{i_1,1},\tilde{A}_{i_2,2},\tilde{A}_{i_3,3}\right)$
for distinct $i_1,i_2,i_3$.  Triangles that come from this triple will
receive the label $(i_1,i_2)$.

The method is very similar to that in Section~\ref{sGAMMACASE}.  The
sets $A_{i,j}'$, $A_{i,j}''$, $\tilde{A}_{i,j}$ and $C_i'$ will be
created similarly to before.  But this case is easier not only because
each $c\in C_i$ can be added to any set $A_{i,j}'$ for any
$j\in\{1,2,3\}$ and it is easy to find the parity triangles.  Parity
triangles can be found in uncolored vertices of
$\left(\tilde{A}_{1,1},\tilde{A}_{2,2},\tilde{A}_{3,3}\right)$,
$\left(\tilde{A}_{2,1},\tilde{A}_{3,2},\tilde{A}_{1,3}\right)$ and
$\left(\tilde{A}_{3,1},\tilde{A}_{1,2},\tilde{A}_{2,3}\right)$.
Partitioning the $\tilde{A}_{i,j}$ sets in half uniformly at random,
exchanging the colored vertices and applying
Proposition~\ref{pEASYCASE} finishes the proof. $\EOP$

\subsection{Proof of Proposition~\ref{pSIZES}}
Let $X$ be the set of vertices in $A_1$ that are adjacent to less than
$(1+\theta)t$ vertices in $B_2$.  Computing the densities,
\[ \Delta|A_1||A_2|\geq e\left(A_1,A_2\right)\geq 2t|A_1|-|A_1||B_2|
   +|X|\left[|B_2|-(1+\theta)t\right] . \]
So, it must be true that
\[ |X|\leq|A_1|\frac{|B_2|-2t+\Delta|A_2|}{|B_2|-(1+\theta)t}
   \leq\frac{\Delta}{1-\theta}t . \] 
As a result, $\left|A_i\setminus A_i'\right|\leq
\frac{2\Delta}{1-\theta}t$.
Similarly, $\left|B_i\setminus B_i'\right|\leq
\frac{4\Delta}{1-\theta}t$.  With
$\Delta_1\geq\frac{4\Delta}{1-\theta}$ the proposition is
proven. $\EOP$

\subsection{Proof of Lemma~\ref{lKEY}}
Again, there are a sequence of constants:
\[ \Delta_1'''\ll\delta_1\ll\delta_2\ll\delta_3\ll\delta_4\ll\delta_5
   \ll\delta_6\ll\delta_7\ll\delta_8\ll
   \Delta_2' . \] 
Begin by defining
\[ B_i'=\left\{v\in V_i : \deg_{A_j'}(v)\geq (1/2)(1+\theta)t',
   \;\forall j\neq i\right\},\quad\mbox{for $i=1,2,3$.} \]
Define $C_i'=V_i\setminus (A_i'\cup B_i')$.  Again, using
Proposition~\ref{pSIZES}, we see that $|C_i'|\leq\delta_1 t'$.
Now, we find $3t'-|A_1'|-|A_2'|-|A_3'|$ disjoint triangles in
$(B_1'\cup C_1',B_2'\cup C_2',B_3'\cup C_3')$.  If this is not
possible, Proposition~\ref{pTHETA} gives that $G$ must be
$\Delta_2'$-approximately $\Theta_{3\times 3}(t')$.

If such disjoint triangles exist, then remove them from the graph to
create $B_i''$ and $C_i''$ for $i=1,2,3$.  All that remains to prove
is that there exists a matching, $M$, in $(B_1''\cup C_1'',B_2''\cup
C_2'',B_3''\cup C_3'')$ such that for any triple $\{i_1,i_2,i_3\}$,
there is a matching in $(B_{i_1}''\cup C_{i_1}'',B_{i_2}''\cup
C_{i_2}'')$ of size $|A_{i_3}'|$ with each $c\in C_{i_1}''$ is
adjacent to at least $(1-\theta-\delta_1)t'$ vertices in $A_{i_3}'$.
What we will do is first form triangles that involve the $c$ vertices
and then, because $\delta_1\ll\theta-3/4$, we can see that each
remaining vertex in, say $A_3'$, is adjacent to at least half of the
edges in the portion of the matching that is in $(B_1''\cup C_1'',
B_2''\cup C_2'')$ and each edge of this portion of the matching is
adjacent to at least half of the remaining vertices in $A_3'$.
K\"onig-Hall gives that there must be a covering by triangles.

In order to find this matching, we will randomly partition the sets
$B_i''\cup C_i''$.  Let $B_i''\cup C_i''=S_i(j)\cup S_i(k)$, where
$\{j,k\}=\{1,2,3\}\setminus\{i\}$ and $|S_i(j)|=|A_j'|$ for all
distinct $i$ and $j$.  It is important to take note that with
probability $1-o(1)$, and for all vertices $v$ in the graph,
\[ \deg_{S_i(j)}(v)-\left(\frac{|A_j'|}{|B_i''\cup C_i''|}\right)
   \deg_{B_i''\cup C_i''}(v)\in\left(-o(t'),+o(t')\right) \; . \] 
This is a result of Stirling's inequality~\cite{B1} (see
Section~\ref{sSTIRLING}).

Once the ``$S$'' sets are randomly chosen, it may be necessary to move
the ``$C''$'' vertices.  Let us suppose that $c\in C_i''\cap S_i(j)$
is not adjacent to at least $(1-\theta-\delta_1)t'$ vertices in
$A_j'$.  Then, we will exchange $c$ with a vertex in $B_i''\cap
S_i(k)$ (where $k=\{1,2,3\}\setminus\{i,j\}$).  Do this for all $i$
and all $c\in C_i''$ and then match each moved vertex in $S_i(j)$ with
an arbitrary neighbor in $S_k(j)$ (for $j\neq k$).  Color these edges
red.  There are at most $\delta_2 t'$ red edges in any pair
$(S_i(j),S_k(j))$.  Then, finish by finding a matching between the
uncolored vertices of $(S_i(j),S_k(j))$.  If this is not possible,
then Proposition~\ref{pHALL}, a simple consequence of K\"onig-Hall,
gives that edges can be removed so that the minimum degree condition
holds, but the pairs must be $\delta_3$-approximately $\Theta_{2\times
2}(|A_j'|/2)$.

\begin{proposition}
Let $\epsilon\ll\Delta$ and $G=(V_1,V_2;E)$ be a balanced
bipartite graph on $2M$ vertices such that each vertex is adjacent to
at least $\left(\frac{1}{2}-\epsilon\right)M$ vertices in the
other class.  If $G$ has no perfect matching, then some edges can be
deleted so that the minimum degree condition is maintained and $G$ is
$\Delta$-approximately $\Theta_{2\times 2}(M/2)$. 
\label{pHALL}
\end{proposition}

If, with probability at least $2/3$, the pair $(S_i(j),S_k(j))$ has
such a matching, then we complete the triangle cover via K\"onig-Hall
and the proof is complete.  Otherwise, with probability at least
$1/3$, the pair $(S_i(j),S_k(j))$ is $\delta_3$-approximately
$\Theta_{2\times 2}(|A_j'|/2)$.  Since this is true and
$\delta_3\ll\delta_4\ll\delta_5$, $(B_i''\cup C_i'',B_k''\cup C_k'')$
itself is $(\delta_4,\delta_5)$-approximately $\Theta_{2\times 2}(t')$.

We want to show that, unless all three pairs are
$(\delta_6,\delta_7)$-approximately $\Theta_{2\times 2}(t')$, the
matching $M$ exists.  Without loss of generality, suppose that
$(B_2''\cup C_2'',B_3''\cup C_3'')$ is not
$(\delta_6,\delta_7)$-approximately $\Theta_{2\times 2}(t')$.  Then
choose ``$S$'' sets as before and move the ``$C$'' vertices as before.
A matching exists among the uncolored vertices of $(S_1(3),S_2(3))$
that involves all but $O(\delta_5)t'$ vertices.  But then exchange
vertices -- outside of this matching -- in $S_2(3)$ with vertices in
$S_2(1)$ so that $M$ can be completed.  If necessary, do the same with
$S_3(2)$ and $S_3(1)$.  Color the edges formed by the switching red.
If there does not exist a matching among the uncolored vertices in
$(S_2(1),S_3(1))$, then, as before, we must have that $(B_2''\cup
C_2'',B_3''\cup C_3'')$ is $(\delta_6,\delta_7)$-approximately
$\Theta_{2\times 2}(t')$, a contradiction.  So, each pair $(B_i''\cup
C_i'',B_k''\cup C_k'')$ must be $(\delta_6,\delta_7)$-approximately
$\Theta_{2\times 2}(t')$.  

The objective is to show that the subsets of vertices involved in
forming the $(\delta_6,\delta_7)$-approximately $\Theta_{2\times
2}(t')$ must have a trivial intersection.  Write $(B_i''\cup
C_i'',B_j''\cup C_j'')$ as $\left(P_{i\rto j}(a)\cup P_{i\rto
j}(b),(P_{j\rto i}(a)\cup P_{j\rto i}(b)\right)$ where each of the
``$P$'' sets are of size $\delta_6$-approximately $t'$ and
\[ d\left(P_{i\rto j}(a),P_{j\rto i}(b)\right),
   d\left(P_{i\rto j}(b),P_{j\rto i}(a)\right))<\delta_7 . \]
Suppose, without loss of generality, that 
\begin{eqnarray*}
   \left|P_{3\rto 1}(a)\cap P_{3\rto 2}(a)\right|, 
   \left|P_{3\rto 1}(a)\cap P_{3\rto 2}(b)\right|, & & \\ 
   \left|P_{3\rto 1}(b)\cap P_{3\rto 2}(a)\right|, 
   \left|P_{3\rto 1}(b)\cap P_{3\rto 2}(b)\right| & \geq & \delta_8 t' .
\end{eqnarray*} 
Then, as in the paragraph above, we can simply choose ``$S$'' sets of
appropriate size at random.  Exchange vertices so as to force a
matching in $(S_1(3),S_2(3))$ and then, since the intersections of the
``$P$'' sets are so large, it is easy to exchange vertices in
$B_3''\cup C_3''$ so that matchings are forced in both
$(S_1(2),S_3(2))$ and $(S_2(1),S_3(1))$.  Using K\"onig-Hall to
complete the covering by triangles gives us a contradiction.

Since, within a tolerance of $\delta_8 t'$, the ``$P$'' sets coincide,
we may assume that $P_{2\rto 1}(a)$ and $P_{2\rto 3}(a)$ coincide and
that $P_{3\rto 2}(a)$ and $P_{3\rto 1}(a)$ coincide.  Therefore, the
issue is whether $P_{1\rto 2}(a)$ and $P_{1\rto 3}(a)$ coincide or
whether they are virtually disjoint.  If they coincide, then $G$ is
$\Delta_2'$-approximately $\Gamma_3(t')$.  If they are disjoint, then
$G$ is $\Delta_2'$-approximately $\Theta_{3\times 3}(t')$. $\EOP$

\subsection{Stirling's Inequality}
\label{sSTIRLING}
Stirling's inequality (see, for example,~\cite{B1}) is a well-known
result that gives
\[
   b(n,k)\exp\left[-\frac{1}{12}\left(\frac{1}{n-k}
   +\frac{1}{k}\right)\right]\leq {n\choose k}\leq b(n,k)
   \exp\left[\frac{1}{12}\left(\frac{1}{n}\right)\right] .
\]
if $b(n,k)=\frac{n^n}{(n-k)^{(n-k)}k^k} \sqrt{\frac{n}{2\pi k(n-k)}}$.
The proofs of Theorem~\ref{tTRIEXTREME} will use the following
corollary:
\begin{corollary}
If $G$ is a graph on $n$ vertices and $X$ is a set on $\Omega(n)$
vertices then, with $\epsilon\ll p$ and $n$ large enough, if $X'$ is
chosen uniformly from ${X\choose p|X|}$,
\[ \Pr\left\{\left|\deg_{X'}(v)-p\deg_X(v)\right|\leq\epsilon
   n, \forall v\in V(G)\setminus X\right\}\rto 1 \] 
as $n\rto\infty$.
\label{cSTIRLING}
\end{corollary}

\section{$N$ is Not a Multiple of 3}
\label{sNOTDIVIS}
We have proven the theorem for the case where $N/3$ is an integer.
The other cases come as a corollary.

Let $t$ be an integer so that $N=3t+1$ and let $N_0=3t$ be large
enough so that Theorem~\ref{tTRI} is true for all multiples of 3
larger than $N_0$.  Remove any triangle from $G$ to form the graph
$G'$.  Then, since every vertex in $G$ is adjacent to at least
$\left\lceil 2t+2/3\right\rceil = 2t+1$ vertices in each of the other
classes of $G$, every vertex in $G'$ is adjacent to at least $2t$
vertices in the other classes of $G'$.  If $G'$ can be covered with
triangles, then clearly $G$ can also.  If $G'=\Gamma_3(t)$, for $t$
odd, then each vertex in $G'$ must be adjacent to both the vertices in
the other vertex classes of $G\setminus G'$.  Therefore, by removing a
triangle in $(A_{1,3},A_{2,3},A_{3,3})$ (with vertex clusters of $G'$
labeled similarly to the diagram in Figure~\ref{efigureTHETA}) and
removing triangles formed by the vertices of $V(G)\setminus V(G')$ and
edges that span the remaining $A$ sets, the resulting graph, $G''$ is
$\Gamma_3(t-1)$, which can be covered with triangles, by the earlier
proof.


The case for $N=3t+2$ is similar.  Remove 2 disjoint triangles.  The
resulting graph can either be covered by disjoint triangles or it is
$\Gamma_3(t)$, for $t$ odd.  In that case we do the same as in the
previous paragraph, forming 5 disjoint triangles and what remains is
the graph $\Gamma_3(t-1)$.

\section{An Open Problem}

An interesting question is that of how to eliminate the phrase ``if
$N$ is large enough.''  It may not be necessary to write another
proof.  In fact, Conjecture~\ref{conjSMALL} would give us a proof of 
Theorem~\ref{tTRI} with $N_0=1$.
\begin{conjecture}
Let $G$ be a graph and $t$ be a positive integer so that both of the
blow-up graphs $G(t)$ and $G(t+1)$ can be covered with triangles.
Then $G$ itself can be covered with triangles.
\label{conjSMALL}
\end{conjecture}

To see that this implies $N_0=1$, suppose that
Conjecture~\ref{conjSMALL} is true.  Furthermore, suppose there is a
balanced tripartite $G$ on $3N$ vertices with the minimum degree
condition, but $G$ cannot be covered with triangles, and
$G\neq\Gamma_3(N/3)$ for $N/3$ odd.  We know by Theorem~\ref{tTRI}
that, for $t\geq N/N_0$, both $G(t)$ and $G(t+1)$ can be covered with
triangles.  This would contradict Conjecture~\ref{conjSMALL}.

\bibliography{magyar}
\bibliographystyle{plain}

\end{document}